\documentclass{article}
\usepackage{amssymb}
\usepackage{amsmath}
\usepackage{amsfonts}

\newtheorem{theorem}{Theorem}

\newtheorem{definition}[theorem]{Definition}

\newenvironment{proof}[1][Proof]{\noindent\textbf{#1.} }{\ \rule{0.5em}{0.5em}}
\input{tcilatex}
\begin{document}

\title{\textbf{CHARACTERIZATIONS OF SPECIAL CURVES}}
\author{\textbf{Yusuf YAYLI}\thanks{%
Ankara University, Faculty of Science, Department of Mathematics, Ankara,
TURKEY}\textbf{, Semra SARACOGLU}\thanks{%
Siirt University, Faculty of Science and Arts, Department of Mathematics,
Siirt, TURKEY}}
\maketitle

\begin{abstract}
In this study, the new characterizations of special curves are investigated
without using the curvatures of these special curves: general helices, slant
helices, Bertrand curves, Mannheim curves. The curvatures are given by the
help of the norms of the derivatives of Frenet vectors.

\textbf{AMS Subj. Class.:} 53A04, 53A05.

\end{abstract}

\section{Introduction}

As is well known fundamental structure of differential geometry is the
curves. Within the process, most of classical differential geometry topics
have been extended to space curves. There are many studies which implies \
different characterizations of these curves. Kula and Yayl\i\ [3] have
investigated spherical images the tangent indicatrix, binormal indicatrix of
a slant helix and obtained that the spherical images are spherical helices.
By defining slant helices and conical geodesic curves, Izumiya and Takeuchi
[1] have considered geometric inviriants of space curves.

In this study, using some approaches in [1] and [3], we give the new
characterizations of special curves. In these characterizations, the
curvatures of these special curves are not used. We show that the curvatures
can be given by the help of the norms of the derivatives of Frenet vectors.
In this connection, some different theorems are presented.

\section{Preliminaries}

Now, we recall some basic concepts of the differential geometry of curves:

\begin{definition}
We assume that the curve $\alpha $ is parametrized by arclength. Then, $%
\alpha ^{\prime }(s)$ is the unit tangent vector to the curve, which we
denote by $T(s)$ Since t has constant length, $T^{\prime }(s)$ will be
orthogonal to $T(s).$ If $T^{\prime }(s)\neq 0$ then we define principal
normal vector 
\begin{equation}
N(s)=\dfrac{T^{\prime }(s)}{\left\Vert T^{\prime }(s)\right\Vert }
\end{equation}%
and the curvature%
\begin{equation*}
\kappa (s)=\left\Vert T^{\prime }(s)\right\Vert
\end{equation*}%
So far, we have
\end{definition}

\begin{equation}
T^{\prime }(s)=\kappa (s).N(s)
\end{equation}%
If $\kappa (s)=0,$ the principal normal vector is not defined. If $\kappa
(s)\neq 0$ then the binormal vector b(s) is given by%
\begin{equation*}
B(s)=T(s)\times N(s)
\end{equation*}%
Then $\left\{ T(s),N(s),B(s)\right\} $ form a right -handed orthonormal
basis for IR$^{3}$. In summary, for the derivatives of Frenet frame, the
Frenet-Serret formulae can be given as [5]:%
\begin{eqnarray}
T^{\prime }(s) &=&\kappa (s).N(s) \\
N^{\prime }(s) &=&-\kappa (s)T(s)+\tau (s)B(s) \\
B^{\prime }(s) &=&-\tau (s)N(s)
\end{eqnarray}%
Here we denote the curvature of the curve $\alpha $ by $\kappa (s)$ and the
torsion of the curve $\alpha $ by $\tau (s).$

\begin{definition}
Let $\alpha $ be a unit speed regular curve in Euclidean 3-space with Frenet
vectors $T,$ $N$ and $B.$ The unit tangent vectors along the curve $\alpha $
generate a curve $\widetilde{T}$ on the sphere of radius 1 about the origin.
The curve $\widetilde{T}$ is called the sphetical indicatrix of $T$ or more
commonly, $\widetilde{T}$ is called the tangent indicatrix of the curve $%
\alpha .$ If $\alpha =\alpha (s)$ is a natural representation of $\alpha ,$
then $\widetilde{T}(s)=T(s)$ will be a representation of $\widetilde{T}$.
Smiliarly one considers the principal normal indicatrix $\widetilde{N}=N(s)$
and binormal indicatrix $\widetilde{B}=B(s)$ [3, 6]

\begin{theorem}
The curve $\alpha $ is a general helix if and only if $\dfrac{\tau }{\kappa }%
(s)=$ constant. If $\kappa (s)\neq 0$ and $\tau (s)$ are constant, it is
called as circular helix.
\end{theorem}

\begin{theorem}
Let $\alpha $ be a unit speed space curve with $\kappa (s)\neq 0.$ Then $%
\alpha $ is a slant helix if and only if 
\begin{equation*}
\sigma (s)=(\frac{\kappa ^{2}}{(\kappa ^{2}+\tau ^{2})^{3/2}}(\dfrac{\tau }{%
\kappa })^{\prime })(s)
\end{equation*}%
is a constant function [1, 3, 4].
\end{theorem}

\begin{theorem}
For a curve $\alpha $ in $E^{3},$ there is a curve $\alpha ^{\ast }$ so that 
$(\alpha ,\alpha ^{\ast })$ is a Mannheim pair [2].
\end{theorem}

\begin{theorem}
Let $\left\{ \alpha ,\alpha ^{\ast }\right\} $ be a Mannheim pair in $E^{3}.$%
The torsion of the curve $\alpha ^{\ast }$ is $\tau ^{\ast }=\dfrac{\kappa }{%
\lambda \tau }.$ Here, $\lambda $ is the distance between corresponding
points of the Mannheim partner curves [2].
\end{theorem}

\begin{theorem}
The curve $\alpha $ is Bertrand curve if and only if $\lambda \kappa +\eta
\tau =1$
\end{theorem}
\end{definition}

\section{CURVES AND THEIR CHARACTERIZATIONS}

Let

\begin{equation}
\alpha :I\rightarrow E^{3}
\end{equation}%
\begin{equation*}
\text{\ \ \ \ \ \ \ }s\mapsto \alpha (s)
\end{equation*}%
be unit speed curve with Frenet vectors $T,N,B$ and with non-zero curvatures 
$\kappa $ and $\tau $ in $%
\mathbb{R}
^{3}.$

In this section, using tangent indicatrix, principal normal indicatrix and
binormal indicatrix of the curve $\alpha ,$ some characterizations have been
given as follows:

\begin{theorem}
The curve $\alpha $ is general helix if and only if $\dfrac{\left\Vert
B^{\prime }\right\Vert }{\left\Vert T^{\prime }\right\Vert }$ is constant.

\begin{proof}
It can be easily seen that

if $T^{\prime }=\kappa N$ then $\left\Vert T^{\prime }\right\Vert =$ $\kappa 
$ and if $B^{\prime }=\tau N$ then $\left\Vert B^{\prime }\right\Vert =\tau
. $ The ratio%
\begin{equation*}
\dfrac{\left\Vert B^{\prime }\right\Vert }{\left\Vert T^{\prime }\right\Vert 
}=\dfrac{\tau }{\kappa }
\end{equation*}

is constant. This completes the proof.

\begin{theorem}
Let the Frenet frame of the spherical tangent indicatrix $\widetilde{T}$ of
the curve $\alpha $ be $\left\{ T,N,B\right\} .$ The curve $\alpha $ is
slant helix if and only if 
\begin{equation*}
\dfrac{\left\Vert D_{T}B\right\Vert }{\left\Vert D_{T}T\right\Vert }%
=cons\tan t
\end{equation*}
\end{theorem}

\begin{theorem}
Let the Frenet frame of the spherical binormal indicatrix $\widetilde{B}$ of
the curve $\alpha $ be $\left\{ \overset{\ast }{T},\overset{\ast }{N},%
\overset{\ast }{B}\right\} .$ The curve $\alpha $ is slant helix if and only
if 
\begin{equation*}
\dfrac{\left\Vert D_{T}\overset{\ast }{B}\right\Vert }{\left\Vert D_{T}%
\overset{\ast }{T}\right\Vert }=cons\tan t
\end{equation*}
\end{theorem}
\end{proof}
\end{theorem}

\begin{theorem}
The curve $\alpha $ is Bertrand curve if and only if 
\begin{equation*}
\lambda \left\Vert T^{\prime }\right\Vert +\eta \left\Vert B^{\prime
}\right\Vert =1
\end{equation*}
\end{theorem}

\begin{theorem}
The curve $\alpha $ is Mannheim curve if and only if 
\begin{eqnarray*}
\dfrac{\left\Vert T^{\prime }\right\Vert }{\left\Vert W\right\Vert ^{2}} &=&%
\dfrac{\left\Vert T^{\prime }\right\Vert }{\left\Vert N^{\prime }\right\Vert
^{2}}=\lambda \\
&=&cons\tan t
\end{eqnarray*}

where the Darboux vector is $W=\tau T+\kappa B.$
\end{theorem}

\begin{theorem}
Let $\alpha $ be a geodesic curve on the surface $M.$ The curve $\alpha $ is
helix on $E^{3\text{ }}$if and only if 
\begin{equation*}
\dfrac{\left\Vert D_{T}Y\right\Vert }{\left\Vert D_{T}T\right\Vert }%
=cons\tan t
\end{equation*}

\begin{proof}
Let $\kappa _{g},\kappa _{n},t_{r}$ be the geodesic curvature, asymptotic
and curvature line respectively. Here it can be easily given that 
\begin{eqnarray*}
\left\Vert D_{T}T\right\Vert &=&\kappa _{n}^{2}+\kappa _{g}^{2} \\
\left\Vert D_{T}Y\right\Vert &=&\kappa _{n}^{2}+t_{r}^{2}
\end{eqnarray*}%
Thus,%
\begin{equation*}
\dfrac{\left\Vert D_{T}Y\right\Vert }{\left\Vert D_{T}T\right\Vert }=\frac{%
1+\left( \dfrac{t_{r}}{\kappa _{n}}\right) ^{2}}{1+\left( \dfrac{\kappa _{g}%
}{\kappa _{n}}\right) ^{2}}
\end{equation*}%
If the curve $\alpha $ is geodesic then $\kappa _{g}=0.$ In that case, 
\begin{equation*}
\dfrac{\left\Vert D_{T}Y\right\Vert }{\left\Vert D_{T}T\right\Vert }%
=cons\tan t
\end{equation*}
\end{proof}
\end{theorem}

\begin{theorem}
Let $\alpha $ be a asymptotic curve on the surface $M.$The curve $\alpha $
is helix on $E^{3\text{ }}$if and only if

\begin{equation*}
\dfrac{\left\Vert D_{T}N\right\Vert }{\left\Vert D_{T}T\right\Vert }%
=cons\tan t
\end{equation*}

\begin{proof}
Smiliarly in Proof 13, it can be easily seen that 
\begin{eqnarray*}
\left\Vert D_{T}T\right\Vert &=&\kappa _{n}^{2}+\kappa _{g}^{2} \\
\left\Vert D_{T}N\right\Vert &=&\kappa _{g}^{2}+t_{r}^{2}
\end{eqnarray*}%
In that case,%
\begin{eqnarray*}
\dfrac{\left\Vert D_{T}N\right\Vert }{\left\Vert D_{T}T\right\Vert } &=&%
\frac{\kappa _{g}^{2}+t_{r}^{2}}{\kappa _{n}^{2}+\kappa _{g}^{2}} \\
&=&\frac{1+\left( \dfrac{t_{r}}{\kappa _{g}}\right) ^{2}}{1+\left( \dfrac{%
\kappa _{n}}{\kappa _{g}}\right) ^{2}}
\end{eqnarray*}%
If the curve $\alpha $ is asymptotic, then $\kappa _{n}=0.$ Thus,%
\begin{equation*}
\dfrac{\left\Vert D_{T}N\right\Vert }{\left\Vert D_{T}T\right\Vert }%
=cons\tan t
\end{equation*}
\end{proof}
\end{theorem}

\begin{theorem}
Let $\left\{ \alpha ,\alpha ^{\ast }\right\} $ be a Mannheim pair. The
torsion of the curve $\alpha ^{\ast }$ is 
\begin{equation*}
\tau ^{\ast }=\dfrac{\left\Vert D_{T}\overset{\ast }{T}\right\Vert }{\lambda
\left\Vert D_{T}\overset{\ast }{B}\right\Vert }=\dfrac{\kappa }{\lambda \tau 
}
\end{equation*}
\end{theorem}

\begin{theorem}
Let $\left\{ \alpha ,\alpha ^{\ast }\right\} $ be a Mannheim pair. The curve 
$\alpha ^{\ast }$ is anti-Salkowski curve if and only if $\alpha $ is a
general helix.

\begin{proof}
As it is seen from the equation%
\begin{equation*}
\tau ^{\ast }=\dfrac{\kappa }{\lambda \tau },
\end{equation*}%
if the curve $\alpha ^{\ast }$ is anti-Salkowski curve then $\tau ^{\ast }$
is constant. Then $\dfrac{\kappa }{\tau }$ is constant. If $\dfrac{\kappa }{%
\tau }$ is constant then the curve $\alpha $ is general helix. This
completes the proof.

\begin{theorem}
The axis of the accompanying screw-motion at a point $c(0)$ is the line in
the direction of the Darboux vector 
\begin{equation*}
W(0)=\tau (0)T(0)+\kappa (0)B(0)
\end{equation*}%
through the point 
\begin{equation*}
P(0)=c(0)+\frac{\kappa (0)}{\kappa ^{2}(0)+\tau ^{2}(0)}N(0)
\end{equation*}%
It can be shown that under these circumtances the tangent to the curve which
passes through all of these points, namely 
\begin{equation*}
P(s)=c(s)+\frac{\kappa (s)}{\kappa ^{2}(s)+\tau ^{2}(s)}N(s)
\end{equation*}%
is proportional to $W(s)$ if and only if $\dfrac{\kappa }{\kappa ^{2}+\tau
^{2}}$ is constant [8].

\begin{proof}
It can be easily seen that%
\begin{equation*}
P^{\prime }(s)=\frac{\tau ^{2}}{\kappa ^{2}+\tau ^{2}}T+\frac{\kappa \tau }{%
\kappa ^{2}+\tau ^{2}}B+(\frac{\kappa }{\kappa ^{2}+\tau ^{2}})^{\prime }N
\end{equation*}%
and 
\begin{equation*}
W(s)=\tau T+\kappa B
\end{equation*}%
Here, if $\dfrac{\kappa }{\kappa ^{2}+\tau ^{2}}$ is constant then 
\begin{eqnarray*}
P^{\prime }(s) &=&\frac{\tau ^{2}}{\kappa ^{2}+\tau ^{2}}T+\frac{\kappa \tau 
}{\kappa ^{2}+\tau ^{2}}B \\
&=&\frac{\tau }{\kappa ^{2}+\tau ^{2}}(\tau T+\kappa B) \\
&=&\frac{\tau }{\kappa ^{2}+\tau ^{2}}W(s)
\end{eqnarray*}%
Thus,%
\begin{equation*}
P^{\prime }(s)=\lambda (s)W(s)
\end{equation*}%
\begin{equation*}
\FRAME{itbpF}{3.6365in}{2.8945in}{0.5725in}{}{}{adsýz.jpg}{\special{language
"Scientific Word";type "GRAPHIC";display "USEDEF";valid_file "F";width
3.6365in;height 2.8945in;depth 0.5725in;original-width
6.6253in;original-height 3.9998in;cropleft "0";croptop "1";cropright
"1";cropbottom "0";filename 'adsýz.JPG';file-properties "XNPEU";}}
\end{equation*}%
\begin{equation*}
Figure.1.\text{The }\tan \text{gent to the curve}
\end{equation*}%
\textbf{Result. }Under these circumtances the tangent to the curve which
passes through all of these points, namely 
\begin{equation*}
P(s)=c(s)+\frac{\kappa (s)}{\kappa ^{2}(s)+\tau ^{2}(s)}N(s)
\end{equation*}%
is proportional to $W(s)$ if and only if $c(s)$ is a Mannheim curve.
\end{proof}
\end{theorem}
\end{proof}
\end{theorem}

\section{CONCLUSIONS}

The starting point of this study is to develope some important
characterizations of special curves by using the curvatures that given by
the help of the norms of the derivatives of Frenet vectors. At this time, it
is obtained that the tangent to the curve is proportional to the Darboux
vector if and only if the curve is Mannheim curve. Additionally, different
theorems have showed that there is a relation between the curvatures and the
norms of the derivatives of Frenet vectors.

We hope that this study will gain different interpretation to the other
studies in this field.

\textbf{References}

\begin{itemize}
\item[{\textbf{[1]}}] Izumiya, S., Takeuchi, N. New Special Curves and
Developable Surfaces, Turk J Math 28(2004), 153-163.

\item[{[\textbf{2}]}] Orbay, K., Kasap, E. On Mannnheim Partner Curves in $%
E^{3}$, International Journal of Physical Sciences, Vol.4(5), pp.261-264,
May, 2009.

\item[{\textbf{[3]}}] Kula, L., Yayl\i , Y. On Slant Helix and its Spherical
Indicatrix, Applied Mathematics and Computation 169(2005) 600-607.

\item[{\textbf{[4]}}] Kula, L., Ekmek\c{c}i, N., Yayl\i , Y. and \.{I}%
larslan, K. Characterizations of Slant Helices in Euclidean 3-space, Turk J
Math 34(2010), 261-273.

\item[{\textbf{[5]}}] Shifrin, T.,.Differential Geometry: A First Course in
Curves and Surfaces (Preliminary Version), University of Georgia, 2010.

\item[{\textbf{[6]}}] Struik, D. J.: Lectures on Classical Differential
Geometry, Dover, New York, 1988.

\item[{\textbf{[7]}}] Wang, F., Liu, H. Mannheim Partner Curves in 3-space,
Proceedings of the Eleventh International Workshop on Differential Geometry
11(2007) 25-31.

\item[{\textbf{[8]}}] Wolfgang K., Differential Geometry
(Curves-Surfaces-Manifolds), American Mathematical Society, 2002.
\end{itemize}

\end{document}